\documentclass[12pt]{article}

\usepackage{amsmath}
\usepackage{amsfonts}
\usepackage{url}
\usepackage{tikz}
\usepackage{endnotes}
\usepackage{framed}
\usepackage{graphicx}
\usepackage{caption}

\usepackage{titlesec}
\titleformat{\section}{\normalsize\bfseries}{\thesection}{1em}{}
\titleformat{\subsection}{\normalsize\bfseries}{\thesubsection}{1em}{}

\begin{document}

\setlength{\parindent}{0in}
\setlength{\parskip}{0.4in plus0.2in minus0.2in}
\setlength{\voffset}{0in}
\setlength{\topmargin}{0in}
\setlength{\headheight}{0in}
\setlength{\headsep}{0in}
\setlength{\footskip}{0.5in}

\noindent

\bibliographystyle{plain}

\begin{center}
{\large\textbf{A Robust Approximation to\\ a Lambert-Type Function}}\\
--------------------------------------------------------------------\\
Ken Roberts\footnote{Physics and Astronomy Department, 
Western University, London, Canada, krobe8@uwo.ca}
\\
April 8, 2015
\end{center}

\begin{abstract}
The function $y = g(x) = \mathrm{log}\big(W(e^x)\big)$,
where $W()$ denotes the Lambert W function,
is the solution to the equation $y + e^y = x$.  
It appears in various problem situations, for instance
the calculation of current-voltage curves for solar cells.
A direct calculation of $g(x)$ may be inaccurate 
because of arithmetic underflow or overflow.
We present a simple algorithm for calculating 
$g(x)$ that is robust, in that it will work for 
almost all $x$ values which are representable
in the arithmetic of one's chosen computer language.
The algorithm does not assume that the chosen
computer language implements the Lambert W function.
\end{abstract}

\section{Introduction}

\begin{figure}
\makebox[\textwidth][c]{\includegraphics[scale=0.40]
   {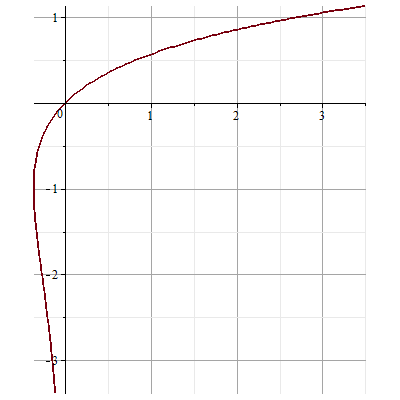}}
\captionsetup{format=hang}
\caption{\label{fig-lamw-graph}Graph of 
$y = Lambert W(x)$ for both $x$ and $y$ real.\\
The graph shows all pairs of real values of $x,y$ \\
which satisfy $y = W(x)$; that is, satisfy $x = y \, e^y$\\
Necessarily $x \ge -1/e$. \\
If $-1/e < x < 0$, then there are two possible $y$ values\\
corresponding to the two branches of Lambert W function.\\
The upper branch, with $y \ge -1$, is the principal branch.\\
If $x$ is real and positive, the only real value of $y = W(x)$ \\
is the positive value of $y$ on the principal branch.
}
\end{figure}

The Lambert W function $w = f(z)$ is the solution of
$w \, e^w = z$, for complex $w$ and $z$.
It can be considered as the multi-branch inverse of the
conformal map $w \to w\,e^w = z$
between the complex $w$-plane and the
complex $z$-plane.
When $w$ and $z$ are restricted to having real values,
the graph of the Lambert W function is as shown in
figure \ref{fig-lamw-graph}. 
For further background regarding the Lambert W function,
see \cite{Corless-1996,Valluri-2000}.

Some problem situations, for instance the modeling of
current and voltage in diodes or solar cells, 
reduce to an implicit equation which can be solved
explicitly by means of the Lambert W function.
As a simple example, consider the implicit equation
\begin{equation}
\label{eqimplicit}
  a \, e^{bU} + bU = V
\end{equation}
where $a, b$ are positive real numbers,
parameters of the model, and $U, V$ are real variables.\footnote{ 
Equation (\ref{eqimplicit}) is a considerable simplification,
in order to have an example in mind.  
The typical solar cell model has four or five parameters.  See \cite{Nelson-2003}, page 14, for instance.}
The corresponding explicit solution for $U$ as a function of $V$ is
\begin{equation}
\label{eqexplicit}
    U = \frac{1}{b} \, \mathrm{log}\Big( \frac{1}{a} \, W(a \, e^V) \Big)
\end{equation}
Here $\mathrm{log}()$ denotes the principal branch
of the natural logarithm function,
and $W()$ denotes the principal branch of the Lambert W function.

A typical task might be, given values of the model parameters, 
to draw a graph of $U$ as a function of $V$.
For that task, it is computationally efficient
to have, instead of the implicit equation (\ref{eqimplicit}),
the explicit solution (\ref{eqexplicit})
for $U$ in terms of $V$ and the model parameters.

Another task might be to estimate the model parameter
values $a, b$ which best fit experimental observations
of $U, V$ pairs.  
For that task, one wants to have an understanding of how
varying the model parameters will affect the $U$-$V$ curves.

Still another task might be to determine the relationships
among the model parameters which correspond to having
an extremum of a function of $U$ or $V$.
This task also requires an understanding of how varying
the model parameters will affect the $U$-$V$ curves,
but it will also be helpful if the formula utilized has partial
derivatives with respect to all the parameters.

The expression in equation (\ref{eqexplicit}) is analytic, 
so is well-behaved
with respect to its argument and the model parameters.
It can be repeatedly differentiated, its extrema lie
at stationary points, and so on.  
Because one is working with real values for $U$ and $V$ and
the positive real model parameters, the argument for the
Lambert W function evaluation in equation (\ref{eqexplicit}) is positive,
so the principal branch of the Lambert W function is being used.
Moreover, one is
using the principal branch of the natural logarithm function.
Everything is single valued in this expression.

However, there can be a numerical difficulty in
performing computations with equation (\ref{eqexplicit}).
Making the substitutions (coordinate changes)
$y = bU+\mathrm{log}(a)$ and $x = V + \mathrm{log}(a)$, 
the computations involve an evaluation of the function
\begin{equation}
\label{eqlwefcn}
  y = g(x) = \mathrm{log}\Big(W(e^x)\Big).
\end{equation}

One difficulty which can arise, depending upon the
computer programming language used, 
is numerical underflow or overflow,
related to the evaluation of an exponential of $x$
where $x$ (negative or positive) has a large magnitude.
The value of $x$ must therefore be restricted to the
set of logarithms of floating point numbers whose exponents can
be accurately represented in the arithmetic facility 
of the computer language.
A second difficulty which can arise is that
the best computer programming language for the
rest of one's problem solution may not have a
built-in Lambert W function evaluator.

The purpose of this note is to address those two
difficulties.  We will describe a simple procedure, 
which can be implemented in any programming
language with floating point arithmetic, for the
robust calculation of the function 
$y = \mathrm{log}(W(\mathrm{exp}(x)))$.
The procedure is valid for essentially any real
value of $x$ which is representable in the
programming language.

\section{Description of the Function 
$y = g(x) = \mathrm{log}(W(\mathrm{exp}(x)))$}

We may consider the function 
$y = g(x) = \mathrm{log}(W(\mathrm{exp}(x)))$
as a transformation of the Lambert W function, with a
change of representation or coordinate space.
For clarity, we will restrict to real arguments.
We can think of $W(u)$,
the principal branch of the Lambert W function, 
as a mapping of the positive real line to itself.
The function $u = h(x) = \mathrm{exp}(x)$ maps the
whole real line to the positive real line,
and its inverse $h^{-1}(u) = \mathrm{log}(u)$,
the principal branch of the natural logarithm,
maps the positive real line to the whole real line.
In this interpretation, the function 
$y = g(x) = \mathrm{log}(W(\mathrm{exp}(x)))$ 
is the composition of functions
$g = h^{-1} \circ W \circ h$, 
and it maps the whole real line to the whole real line.

Suppose that $y = \mathrm{log}(W(\mathrm{exp}(x)))$.
Taking exponentials (ie, applying the function $h$
to both sides of the equation) gives
\begin{equation*}
  e^y = W(e^x).
\end{equation*}
That is, using the definition of the Lambert W function,
\begin{equation*}
  e^y \, e^{e^y} = e^x
\end{equation*}
or
\begin{equation*}
  e^{y + e^y} = e^x.
\end{equation*}
Taking logarithms (ie, applying the function $h^{-1}$
to both sides) gives
\begin{equation}
\label{eqyey}
   y + e^y = x.
\end{equation}

Equation (\ref{eqyey}) is a simple equation structure, 
as simple as the
Lambert W defining equation structure 
\begin{equation}
\label{eqwew}
  w \, e^w = z.
\end{equation}
In fact, equation (\ref{eqyey}) is just equation 
(\ref{eqwew}) in another coordinate system.

When we are evaluating 
$y = g(x) = \mathrm{log}(W(\mathrm{exp}(x)))$
as the solution to equation (\ref{eqyey}),
we are just evaluating the Lambert W function.
There is an important difference, however:
The evaluation of $g(x)$ does not involve much
risk of underflow
or overflow in the numerical representation of
the computer language.

Since $y = g(x)$ satisfies $g(x) + e^{g(x)} = x$,
the first derivative of $g(x)$ satisfies
\begin{equation*}
    g'(x) + e^{g(x)} \, g'(x) = 1
\end{equation*}
or
\begin{equation*}
    g'(x) = \frac{1}{1 + e^{g(x)}} > 0
\end{equation*}

The second derivative of $g(x)$ satisfies
\begin{equation*}
   g''(x) = -\frac{e^{g(x)} \, g'(x)}{(1+e^{g(x)})^2} < 0
\end{equation*}

\begin{figure}
\makebox[\textwidth][c]{
  \includegraphics[scale=0.60]{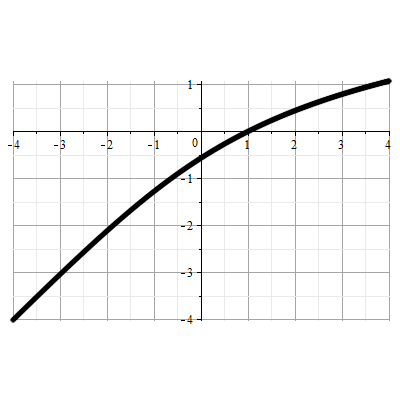}}
{
\captionsetup{format=hang}
  \caption{
\label{fig-logwexp1}Graph of 
  $y = g(x) = \mathrm{log}(W(\mathrm{exp}(x)))$\\
for moderate magnitudes of the argument.
}}
\end{figure}

\begin{figure}
\makebox[\textwidth][c]{
  \includegraphics[scale=0.60]{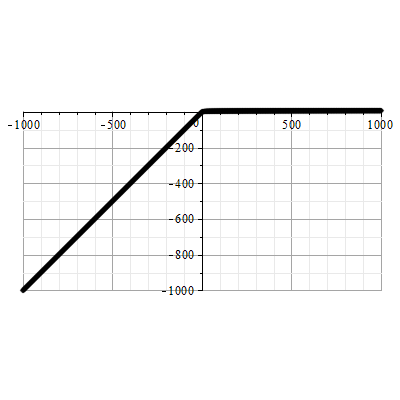}}
{
\captionsetup{format=hang}
  \caption{
\label{fig-logwexp2}Graph of 
  $y = g(x) = \mathrm{log}(W(\mathrm{exp}(x)))$\\
for large magnitudes of the argument.
}}
\end{figure}

Figure \ref{fig-logwexp1} shows the function 
$y = g(x)$ for moderate values of the argument $x$,
that is, for $-4 < x < 4$.
Figure \ref{fig-logwexp2} shows the same function
for larger values of the argument $x$,
that is, for $-1000 < x < 1000$.

One can see from these graphs that the function
$g(x)$ behaves like $x$ when $x$ is much less than
0, and behaves like $\mathrm{log}(x)$ when $x$ is
much more than 0.  
For values of $x$ around 0, there is a smooth blend
between the two behaviors, with $g(1) = 0$.
The function $g(x)$ is strictly monotonic
increasing, as it has a positive first derivative.
It curves downward, as it has a negative second derivative.

One can further see, from figure \ref{fig-logwexp2},
that when the range of the argument $x$ is large,
the graph of $y = g(x)$ looks like it has a sharp corner
at the origin.
Actually, as figure \ref{fig-logwexp1} illustrates,
the graph does not really have a sharp corner.
Nonetheless, at a suitable distance (large scale),
one has in the function $g$ a useful smooth function
for representing a function which has a step in its
derivative.

\section{Algorithm for Calculating 
$y = g(x) = \mathrm{log}(W(\mathrm{exp}(x)))$}

In the terminology of H. Kuki (\cite{Kuki-1966} page 23),
the function $g(x)$ is contracting.  That is,
$|\Delta g(x)| \le |\Delta x|$, or $|g'(x)| \le 1$,
for all $x$ values in its argument domain.
That means the task of finding an estimate for
$g(x)$ given $x$ is relatively stable.
A slight change in $x$ (noise in the input)
will produce only a slight change in $g(x)$.
The only challenges in developing a formula
to estimate $g(x)$, given $x$, are finding
an appropriate algorithm, coding the
sequence of calculations to avoid unwanted
cancellation, and being reasonably efficient
in the number of computations performed.

A suitable algorithm can be an initial estimate,
followed by some number of iterations of a
refinement.  
Halley's method is used to perform refinements
because it has cubic convergence, and the
derivatives involved can be calculated efficiently.

Given any fixed real number $x$, we wish to find a
real number $y$ such that
\begin{equation*}
    h(y) = y + e^y - x
\end{equation*}
is zero.  
The first and second derivatives of $h(y)$ are needed
for Halley's method.  They are
\begin{eqnarray*}
    h'(y) &=& 1 + e^y
\\
    h''(y) &=& e^y
\end{eqnarray*}
and hence are particularly easy to calculate.
Once one has $e^y$ from the calculation of
$h(y)$, the derivatives are also at hand.
It is also necessary, in order to use Halley's
method, that the first derivative is non-zero;
that is the case for $h'(y) = 1 + e^y$.

As an initial estimate $y_0$, we choose to use 
$y_0 = x$ for $x \le -e$,
and $y_0 = \mathrm{log}(x)$ for $x \ge e$.
For $-e < x < e$, we linearly interpolate
between the two values $-e$ and 1.
This is an extraordinarily crude initial estimate,
but it is sufficient, since Halley's method is
very robust and rapidly convergent in this application.

The general iteration formula for Halley's method is
\begin{equation*}
  y_{n+1} = y_n - \frac{2 h(y_n) h'(y_n)}{2 h'(y_n)^2 - h(y_n)h''(y_n)}
\end{equation*}
In this particular case  $h(y) = y + e^y - x$ and the
iteration formula becomes
\begin{equation}
\label{eqiter}
  y_{n+1} = y_n - 
  \frac{2 (y_n + e^{y_n} - x) (1 + e^{y_n})}
         {2 (1 + e^{y_n})^2 - (y_n + e^{y_n} - x) e^{y_n}}.
\end{equation}

The details of coding depend upon the computer language.
It will be efficient to evaluate $e^{y_n}$ only once per iteration.
All other computations are straightforward arithmetic.
When evaluating the denominator of the adjustment in
the iteration equation (\ref{eqiter}), there is little risk
of cancellation resulting from the subtraction, as
the first term in the denominator is larger than the
second term.

In practice, just a few iterations suffice to give a good
result.   For arguments in $-10^6 \le x \le 10^6$,
four iterations of Halley's method 
reduce the absolute error to less
than $10^{-80}$.
The actual coding can use a convergence criterion,
based upon the desired maximum error in the
estimate of function value, to
determine how many iterations to perform.
Alternatively, if the precision is fixed by the computer
language's arithmetic representation or by the needs
of the application situation, then one can determine
how many iterations of the refinement will suffice,
and perform only that number, omitting the final
redundant iteration which verifies the convergence.
This technique, due to H. Kuki as seen in his 
algorithms for computing square root (see \cite{Kuki-1966},
pages 49-50 and 135-136), probably deserves a name.
Perhaps it should be called ``Kuki's convergence non-test".

\section{Acknowledgments}

The author had the good fortune and honor to work for
Hirondo Kuki in 1968-69, and thanks him for his guidance,
support and friendship. 
He also thanks S. R. Valluri for an introduction to the
Lambert W function and the many interesting problems
associated with its properties, and in particular for
a stimulating discussion of the topic of this note.
He thanks Mark Campanelli for suggesting Halley's
method for iterative refinement in solar cell calculations.


\end{document}